\documentclass[11pt]{article}

\usepackage{amssymb,amsmath,amsfonts,amsthm,a4wide}

\newtheorem{theorem}{Theorem}
\newtheorem{lemma}{Lemma}
\newtheorem{assumption}{Assumption}
\newtheorem*{remark}{Remark}
\newtheorem*{remarks}{Remarks}
\newtheorem*{lemma*}{Lemma}

\newcommand{\RR}{\mathbb{R}}
\newcommand{\CC}{\mathbb{C}}

\newcommand{\DD}{\Delta}
\newcommand{\Tm}{\sqrt{-\DD + m^2} \, } 
\newcommand{\T}{\sqrt{-\DD} \,  } 
\newcommand{\ie}{i.\,e.}
\newcommand{\eg}{e.\,g.}

\newcommand{\Lploc}[2]{L^{#1}_{\mathrm{loc}}(\RR^{#2})}

\newcommand{\eps}{\varepsilon}
\newcommand{\Me}{\mathcal{M}_\eps}
\newcommand{\Ds}[1]{\mathcal{D}^{#1}}

\numberwithin{equation}{section}

\begin{document}

\title{{\bf Well-Posedness for Semi-Relativistic Hartree Equations of Critical Type}} 
\author{Enno Lenzmann \\ Department of Mathematics, ETH Z\"urich \\ E-Mail: {\tt lenzmann@math.ethz.ch}}
\date{August 30, 2005}
\maketitle

\begin{abstract}
We prove local and global well-posedness for semi-relativistic, nonlinear Schr\"odinger equations $i \partial_t u = \sqrt{-\Delta + m^2} u + F(u)$ with initial data in $H^s(\mathbb{R}^3)$, $s \geq 1/2$. Here $F(u)$ is a critical Hartree nonlinearity that corresponds to Coulomb or Yukawa type self-interactions. For focusing $F(u)$, which arise in the quantum theory of boson stars, we derive a sufficient condition for global-in-time existence in terms of a solitary wave ground state. Our proof of well-posedness does not rely on Strichartz type estimates, and it enables us to add external potentials of a general class. 
\end{abstract}

%
 
%

\section{Introduction} \label{sec-intro}

In this paper we study the Cauchy problem for nonlinear Schr\"odinger equations with kinetic energy part originating from special relativity. That is, we consider the initial value problem for 
\begin{equation} \label{eq-snls}
 i \partial_t u = \Tm u + F(u),  \qquad  (t,x) \in \RR^{1+3}, 
\end{equation}
where $u(t,x)$ is complex-valued, $m \geq 0$ denotes a given mass parameter, and $F(u)$ is some nonlinearity. Here the operator $\sqrt{-\DD + m^2}$ is defined via its symbol $\sqrt{\xi^2 + m^2}$ in Fourier space.

Such ``semi-relativistic'' equations have (though not Lorentz covariant in general) interesting applications in the quantum theory for large systems of self-interacting, relativistic bosons. Equation (\ref{eq-snls}) arises, for instance, as an effective description of {\em boson stars}, see, \eg, \cite{Elgart+Schlein2005,Lieb+Yau1987}, where $F(u)$ is a focusing Hartree nonlinearity given by 
\begin{equation} \label{eq-F}
F(u) =  \big ( \frac{\lambda}{|x|} \ast |u|^2 ) u,
\end{equation}
with some constant $\lambda < 0$ and $\ast$ as convolution. Motivated by this physical example with focusing self-interaction of Coulomb type, we address the Cauchy problem for equation (\ref{eq-snls}) and a class of Hartree nonlinearities including (\ref{eq-F}). In fact, we shall prove well-posedness for initial data $u(0,x) = u_0(x)$ in $H^{s}=H^{s}(\RR^3)$, $s \geq 1/2$; see Theorems \ref{th-lwp}--\ref{th-reg} below.

Let us briefly point out a decisive feature of the example cited in (\ref{eq-F}) above. Apart from its physical relevance, the nonlinearity given by (\ref{eq-F}) leads to an {\em $L^2$-critical} equation as indicated by the fact that the coupling constant $\lambda$ has to be dimensionless. In consequence of this, $L^2$-smallness of the initial datum enters as a sufficient condition for global-in-time solutions. More precisely, we derive for $u_0 \in H^{s}$, $s \geq 1/2$, the following criterion implying global well-posedness
\begin{equation} \label{eq-cond}
 \int_{\RR^3} | u_0(x) |^2 \, dx < \int_{\RR^3} |Q(x)|^2 \, dx .
\end{equation}
This condition holds irrespectively of the parameter $m \geq 0$ in (\ref{eq-snls}); see Theorem \ref{th-gwp} below. Here $Q \in H^{1/2}$ is a positive solution (ground state) for the nonlinear equation
\begin{equation}
\T Q + \big ( \frac{\lambda}{|x|} \ast |Q|^2 ) Q = - Q,
\end{equation}
which gives rise to solitary wave solutions, $u(t,x) = e^{it} Q(x)$, for (\ref{eq-snls}) with $m=0$. In fact, it can be shown that criterion (\ref{eq-cond}) guaranteeing global-in-time solutions in the focusing case is optimal in the sense that there exist solutions, $u(t)$, with $\| u_0 \|_2^2 > \| Q \|_2^2$, which blow up within finite time; see \cite{Lenzmann2005} for a proof.

Furthermore, criterion (\ref{eq-cond}) can be linked with established results as follows. First, it is reminiscent to a well-known condition derived in \cite{Weinstein1983} for global well-posedness of nonrelativistic Schr\"odinger equations with focusing, local nonlinearity (see also \cite{Nawa+Ozawa1992} for Hartree nonlinearities). Second, criterion (\ref{eq-cond}) is in accordance with a sufficient stability condition proved in \cite{Lieb+Yau1987} for the related time-independent problem (\ie, a static boson star); see \cite{Froehlich+Lenzmann2004} for a more details concerning known results on Hartree equations. 
 
We now give an outline of our methods. The proof of well-posedness presented below does {\em not} rely on Strichartz (\ie, space-time) estimates for the propagator, $e^{-it \Tm}$, but it employs sharp estimates (\eg, Kato's inequality (\ref{ineq-kato}) below) to derive local Lipschitz continuity of $L^2$-critical nonlinearities of Hartree type. Local well-posedness then follows by standard methods for abstract evolution equations. Furthermore, global well-posedness is derived by means of a-priori estimates and conservation of charge and energy whose proof requires a regularization method. 

This paper is organized as follows. 
\begin{itemize}

\item In Section \ref{sec-main} we introduce a class of critical Hartree nonlinearities including (\ref{eq-F}). First, we state Theorems \ref{th-lwp} and \ref{th-gwp} that establish local and global well-posedness in energy space $H^{1/2}$ for this class of nonlinearities. In Theorem \ref{th-reg} we extend these results to $H^s$, for every $s \geq 1/2$. Finally, external potentials are included, \ie, we consider 
\begin{equation} \label{eq-VV}
 i \partial_t u = \big ( \Tm + V ) u + F(u) ,  
\end{equation}
where $V : \RR^3 \rightarrow \RR$ is given. In Theorem \ref{th-V} we state local and global well-posedness for (\ref{eq-VV}) with initial datum $u(0,x)=u_0(x)$ in the appropriate energy space. Assumption \ref{ass-V} imposed below on $V$ is considerably weak and implies that $\Tm + V$ defines a self-adjoint operator via its form sum.

\item The main results (\ie, Theorems \ref{th-lwp}--\ref{th-V}) are proved in Section \ref{sec-proofs}.

\item Appendix \ref{sec-app} contains useful facts about fractional derivatives, a discussion of ground states, and some details of the proofs.  

\end{itemize}

\subsubsection*{Notation} 
Throughout this text, the symbol $\ast$ stands for convolution on $\RR^3$, \ie,
\[ (f \ast g)(x) := \int_{\RR^3} f(x-y) g(y) \, dy,\]
and $L^p(\RR^3)$, with norm $\| \cdot \|_p$ and $1 \leq p \leq \infty$, denotes the usual Lebesgue $L^p$-space of complex-valued functions on $\RR^3$. Moreover, $L^2(\RR^3)$ is associated with the scalar product defined by
\[
\langle u, v \rangle := \int_{\RR^3} \overline{u}(x)  v(x) \, dx .
\] 
For $s \in \RR$ and $1 \leq p \leq \infty$, we introduce fractional Sobolev spaces (see, \eg, \cite{Bergh+Loefstroem1976}) with their corresponding norms according to 
\begin{equation*}
H^{s,p}(\RR^3) := \big \{ u \in \mathcal{S}'(\RR^3) : \| u \|_{H^{s,p}} :=  \| \mathcal{F}^{-1} [ (1+\xi^2)^{s/2} \mathcal{F} u ]  \|_p < \infty \big \} ,
\end{equation*} 
where $\mathcal{F}$ denotes the Fourier transform in $\mathcal{S}'(\RR^3)$ (space of tempered distributions). In our analysis, the Sobolev spaces
\[ H^{s}(\RR^3) := H^{s,2}(\RR^3),  \]
with norms $\| \cdot \|_{H^s} := \| \cdot \|_{H^{s,2}}$, will play an important role.

In addition to the common $L^p$-spaces, we also make use of local $L^p$-space, $L^p_\mathrm{loc}(\RR^3)$, with $1 \leq p \leq \infty$, and weak (or Lorentz) spaces, $L^p_w(\RR^3)$, with $1 < p < \infty$ and corresponding norms given by
\[
\| u \|_{p,w} := \sup_\Omega | \Omega |^{-1/p'} \int_\Omega |u(x)| \, dx,
\]  
where $1/p + 1/p' = 1$ and $\Omega$ denotes an arbitrary measurable set with Lebesgue measure $|\Omega| < \infty$; see, \eg, \cite{Lieb+Loss2001} for this definition of $L^p_w$-norms. Note that $L^p(\RR^3) \subsetneq L^p_w(\RR^3)$, for $1 < p < \infty$.

The symbol $\DD = \sum_{i=1}^3 \partial^2_{x_i}$ stands for the usual Laplacian on $\RR^3$, and $\Tm$ is defined via its symbol $\sqrt{\xi^2 + m^2}$ in Fourier space. Besides the operator $\Tm$, we also employ Riesz and Bessel potentials of order $s \in \RR$, which we denote by $(-\DD)^{s/2}$ and $(1-\DD)^{s/2}$, respectively; see also Appendix \ref{sec-app}. 

Except for theorems and lemmas, we often use the abbreviations $L^p = L^p(\RR^3)$, $L^p_w = L^p_w(\RR^3)$, and $H^s = H^s(\RR^3)$. In what follows, $a \lesssim b$ always denotes an inequality $a \leq cb$, where $c$ is an appropriate positive constant that can depend on fixed parameters.

\section{Main Results} \label{sec-main}

We consider the following initial value problem 
\begin{equation} \label{eq-Hartree}
\left \{ \begin{array}{l}  \displaystyle  i \partial_t u = \Tm u +  \big ( \frac{\lambda e^{-\mu |x|}}{|x|} \ast |u|^2 \big ) u,  \\ 
\displaystyle u(0,x) = u_0(x), \quad u : [0,T) \times \RR^3 \rightarrow \CC ,  \end{array} \right .
\end{equation}
where $m \geq 0$, $\lambda \in \RR$, and $\mu \geq 0$ are given parameters. Note that $|\lambda|$ could be absorbed in the normalization of $u(t,x)$, but we shall keep $\lambda$ explicit in the following; see also \cite{Elgart+Schlein2005} for this convention.    

Our particular choice of the Hartree type nonlinearities in (\ref{eq-Hartree}) is motivated by the fact that (\ref{eq-Hartree}) can be rewritten as the following system of equations 
\begin{equation} \label{eq-CS}
\left \{ \begin{array}{l}  i \partial_t u = \Tm u + \Psi u,  \\
\displaystyle (\mu^2-\DD) \Psi = 4 \pi \lambda | u |^2 , \quad u(0,x) = u_0(x), \end{array} \right .  
\end{equation}
where $\Psi= \Psi(t,x)$ is real-valued and $\Psi(t,x) \rightarrow 0$ as $|x| \rightarrow \infty$. This reformulation stems from the observation that $e^{-\mu |x|}/4 \pi |x|$ is the Green's function of $(\mu^2 - \DD)$ in $\RR^3$; see Appendix \ref{app-calc}. System (\ref{eq-CS}) now reveals the physical intuition behind (\ref{eq-Hartree}), \ie, the function $u(t,x)$ corresponds to a ``positive energy wave'' with instantaneous self-interaction that is either of Coulomb or Yukawa type depending on whether $\mu = 0$ or $\mu > 0$, respectively. To prove well-posedness we shall, however, use formulation (\ref{eq-Hartree}) instead, and we refer to facts from potential theory only when estimating the nonlinearity.

\subsection{Local Well-Posedness}

Let us begin with well-posedness in energy space, \ie, we assume that $u_0 \in H^{1/2}$ holds in (\ref{eq-Hartree}). The following Theorem \ref{th-lwp} establishes local well-posedness in the strong sense, \ie, we have existence and uniqueness of solutions, their continuous dependence on initial data, and the blow-up alternative. The precise statements is as follows.

\begin{theorem} \label{th-lwp}
Let $m \geq 0, \lambda \in  \RR$, and $\mu \geq 0$. Then initial value problem (\ref{eq-Hartree}) is locally well-posed in $H^{1/2}(\RR^3)$. This means that, for every $u_0 \in H^{1/2}(\RR^3)$, there exist a unique solution
\[ u \in C^0 \big( [0,T); H^{1/2}(\RR^3) \big ) \cap C^1 \big ([0,T); H^{-1/2}(\RR^3) \big ), \]
and it depends continuously on $u_0$. Here $T  \in (0, \infty]$ is the maximal time of existence, where we have that either $T = \infty$ or $T < \infty$ and $\lim_{t \uparrow T} \| u(t) \|_{H^{1/2}} = \infty$ holds. 
\end{theorem}

\begin{remark} {\em Continuous dependence means that the map $u_0 \mapsto u \in C^0(I; H^{1/2})$ is continuous for every compact interval $I \subset [0,T)$. }
\end{remark}

\subsection{Global Well-Posedness}

The local-in-time solutions derived in Theorem \ref{th-lwp} extend to all times, by virtue of Theorem \ref{th-gwp} below, provided that either $\lambda \geq 0$ holds (corresponding to a repulsive nonlinearity) or $\lambda < 0$ and the initial datum is sufficiently small in $L^2$. 

\begin{theorem}\label{th-gwp}
The solution of (\ref{eq-Hartree}) derived in Theorem \ref{th-lwp} is global in time, \ie, we have that $T = \infty$ holds, provided that one of the following conditions is met.
\begin{enumerate}
\item[i)] $\lambda \geq 0$. 
\item[ii)] $\lambda < 0$ and $\| u_0 \|^2_2 < \| Q \|_2^2$, where $Q \in H^{1/2}(\RR^3)$ is a strictly positive solution (ground state) of
\begin{equation} \label{eq-Q1} 
\T Q + \big ( \frac{\lambda}{|x|} \ast |Q|^2 \big ) Q = - Q.
\end{equation} 
Moreover, we have the estimate $\| Q \|_2^2 > \frac{4}{\pi  |\lambda|}$.
\end{enumerate}
\end{theorem}

\begin{remarks} {\em 1) Notice that condition ii) implies global well-posedness for (\ref{eq-Hartree}) irrespectively of $m \geq 0$.

2) Due to the scaling behavior of (\ref{eq-Q1}), the function $Q_a(x) = a^{3/2} Q(ax)$, with $a > 0$, yields another ground state with $\| Q_a \|_2 = \| Q \|_2$ that satisfies
\begin{equation}
 \T Q_a + \big ( \frac{\lambda}{|x|} \ast |Q_a|^2) Q_a = -a Q_a. 
\end{equation}
We refer to Appendix \ref{app-ground} for a discussion of $Q \in H^{1/2}$. 

3) Condition ii) resembles a well-known criterion derived in \cite{Weinstein1983} for global-in-time existence for $L^2$-critical nonlinear (nonrelativistic) Schr\"odinger equations. 

4) It is shown in \cite{Lenzmann2005} that criterion (\ref{eq-cond}) for having global-in-time solutions in the focusing case is optimal in the sense that there exist solutions, $u(t)$, with $ \| u_0 \|_2^2 > \| Q \|_2^2$, which blow up within finite time. }
\end{remarks}


\subsection{Higher Regularity}

We now turn to well-posedness of (\ref{eq-Hartree}) in $H^s$, for $s \geq 1/2$, which is settled by the following result.
\begin{theorem}\label{th-reg}
For every $s \geq 1/2$, the conclusions of Theorems \ref{th-lwp} and \ref{th-gwp} hold, where $H^{1/2}(\RR^3)$ and $H^{-1/2}(\RR^3)$ in Theorem \ref{th-lwp} are replaced by $H^s(\RR^3)$ and $H^{s-1}(\RR^3)$, respectively.
\end{theorem}

\begin{remark} {\em For $s=1$, this result is needed in \cite{Elgart+Schlein2005} for a rigorous derivation of (\ref{eq-Hartree}) with Coulomb type self-interaction (\ie, $\mu=0$) from many-body quantum mechanics. }
\end{remark}

\subsection{External Potentials}
  
Now we consider the following extension of (\ref{eq-Hartree}) that arises by adding an external potential: 
\begin{equation} \label{eq-HartreeV}
\left \{ \begin{array}{l}  \displaystyle  i \partial_t u = \big ( \Tm + V \big )  u +   \big ( \frac{ \lambda e^{- \mu |x|}}{|x|} \ast |u|^2 \big ) u,  \\ 
u(0,x) = u_0(x), \quad u : [0,T) \times \RR^3 \rightarrow \CC, \end{array} \right  . 
\end{equation}
where $m \geq 0$, $\lambda \in \RR$, $\mu \geq 0$ are given parameters, and $V : \RR^3 \rightarrow \RR$ denotes a preassigned function that meets the following condition.

\begin{assumption} \label{ass-V}
Suppose that $V = V_+ + V_-$ holds, where $V_+$ and $V_-$ are real-valued, measurable functions with the following properties.
\begin{enumerate}
\item[i)] $V_+ \in \Lploc{1}{3}$ and $V_+ \geq 0$.
\item[ii)] $V_-$ is $\sqrt{-\DD}$-form bounded with relative bound less than $1$, \ie, there exist constants $0 \leq a < 1$ and $0 \leq b < \infty$, such that 
\begin{equation*} 
|\langle u, V_- u \rangle | \leq a \langle u, \T u \rangle + b \langle u, u \rangle
\end{equation*}
holds for all $u \in H^{1/2}(\RR^3)$.
\end{enumerate}
\end{assumption}
We mention that Assumption 1 implies that $\Tm + V$ leads to a self-adjoint operator on $L^2$ via its form sum. Furthermore, the energy space given by
\begin{equation} 
X := \big \{ u \in H^{1/2}(\RR^3) :  \int_{\RR^3} V(x) \, |u(x)|^2 \, dx < \infty \big \}
\end{equation}
is complete with norm $\| \cdot \|_X$, and its dual space is denoted by $X^*$. We refer to Section \ref{sec-proof-V} for more details on $\Tm + V$ and $X$. 
 
After this preparing discussion, the extension of Theorems \ref{th-lwp} and \ref{th-gwp} for the initial value problem (\ref{eq-HartreeV}) can be now stated as follows.

\begin{theorem} \label{th-V} Let $m \geq 0, \lambda \in \RR$, $\mu \geq 0$, and suppose that $V$ satisfies Assumption \ref{ass-V}. Then (\ref{eq-HartreeV}) is locally well-posed in the following sense. For every $u_0 \in X$, there exists a unique solution
\[
u \in C^0([0,T); X) \cap C^1([0,T); X^*),
\]  
and it depends continuously on $u_0$.  Here $T  \in (0, \infty]$ is the maximal time of existence such that either $T = \infty$ or $T < \infty$ and $\lim_{t \uparrow T} \| u(t) \|_{X} = \infty$ holds. Moreover, we have that $T= \infty$ holds, if one of the following conditions is satisfied.
\begin{enumerate}
\item[i)] $\lambda \geq 0$.
\item[ii)] $\lambda < 0$ and $ \| u_0 \|_2^2 < (1-a) \|Q \|_2^2$, where $Q$ is the ground state mentioned in Theorem \ref{th-gwp} and $0 \leq a < 1$ denotes the relative bound introduced in Assumption \ref{ass-V}.
\end{enumerate}
\end{theorem} 

\begin{remarks} {\em 1) To meet Assumption \ref{ass-V} for $V_+$, we can choose, for example, $V_+(x) = |x|^{\beta}$, with $\beta \geq 0$; or even super-polynomial growth such as $V_+(x) = e^{x}$. Note that Assumption \ref{ass-V} for $V_-$ is satisfied (by virtue of Sobolev inequalities), if
\[ |V_-(x)| \leq \frac{c}{|x|^{1-\eps}} + d \]
holds for some $0 < \eps \leq 1$ and constants $0 \leq c,d < \infty$. In fact, we can even admit $\eps = 0$ provided that $c < 2/\pi$ holds, as can be seen from inequality (\ref{ineq-kato}) below.     

2) Since we avoid Strichartz estimates in our well-posedness proof below, we only need that $V_+$ belongs to $L^1_\mathrm{loc}$. In contrast to this, compare, for instance, the conditions on $V$ in \cite{Yajima+Zhang2004} for deriving Strichartz type estimates for $e^{-it(-\DD+V)}$ in order to prove local well-posedness for (nonrelativistic) nonlinear Schr\"odinger equations with external potentials. }
\end{remarks}


\section{Proof of the Main Results} \label{sec-proofs}

In this section we prove Theorems \ref{th-lwp}--\ref{th-V}. Although Theorem \ref{th-V} generalizes Theorems \ref{th-lwp} and \ref{th-gwp}, we postpone the proof of Theorem \ref{th-V} to the final part of this section. 


\subsection{Proof of Theorem \ref{th-lwp} (Local Well-Posedness)}

Let $u_0 \in H^{1/2}$ be fixed. In view of (\ref{eq-Hartree}) we put
\begin{equation} \label{eq-AF}
A := \Tm \quad \mbox{and} \quad F(u) := \big ( \frac{\lambda e^{-\mu |x|}}{|x|} \ast |u|^2 \big ) u ,
\end{equation}  
and we consider the integral equation
\begin{equation} \label{eq-integral}
 u(t) = e^{-itA}  u_0 - i \int_0^t e^{-i(t-\tau)A} F(u(\tau)) \, d \tau. 
\end{equation}
Here $u(t)$ is supposed to belong to the Banach space
\begin{equation}
 Y_T := C^0 \big( [0,T); H^{1/2}(\RR^3) \big ),
\end{equation}
with some $T >0$ and corresponding norm $\| u \|_{Y_T} := \sup_{t \in [0,T)} \| u(t) \|_{H^{1/2}}$. The proof of Theorem \ref{th-lwp} is now organized in two steps as follows.


\subsubsection*{Step 1: Estimating the Nonlinearity}

We show that the nonlinearity $F(u)$ is locally Lipschitz continuous from $H^{1/2}$ into itself. This is main point of our argument for local well-posedness and it reads as follows.

\begin{lemma} \label{lem-lip}
For $\mu \geq 0$, the map $J(u) := \big ( \frac{e^{-\mu |x|}}{|x|} \ast |u|^2 \big ) u$ is locally Lipschitz continuous from $H^{1/2}(\RR^3)$ into itself with
\[ \| J(u) - J(v) \|_{H^{1/2}} \lesssim ( \| u \|^2_{H^{1/2}} + \| v \|^2_{H^{1/2}} ) \| u - v \|_{H^{1/2}}, \]
for all $u, v \in H^{1/2}(\RR^3)$. 
\end{lemma}

\begin{proof}[Proof of Lemma \ref{lem-lip}] We prove the claim for $\mu = 0$ and $\mu > 0$ in a common way, so let $\mu \geq 0$ be fixed. For $s \in \RR$, it is convenient to introduce
\[ \Ds{s} := (\mu^2-\DD)^{s/2} . \]
Note that due to the equivalence
\begin{equation*}
 \| u \|_2 + \| \Ds{1/2} u \|_2 \lesssim \| u \|_{H^{1/2}} \lesssim \| u \|_2 + \| \Ds{1/2} u \|_2 , 
 \end{equation*}
it is sufficient to estimate the quantities
\begin{equation*}
  I:= \| J(u) - J(v) \|_2 \quad \mbox{and} \quad II := \| \Ds{1/2} [ J(u) - J(v) ] \|_2, 
 \end{equation*}
where $I$ is needed only if $\mu = 0$. Using now the identity
\begin{equation*}
 J(u) - J(v) = \frac{1}{2} \Big [  \big ( \frac{e^{-\mu |x|}}{|x|} \ast (|u|^2 - |v|^2) \big )(u+v) + \big ( \frac{e^{-\mu |x|}}{|x|} \ast ( |u|^2 + |v|^2) \big ) (u-v) \Big ] 
\end{equation*}
together with H\"older's inequality (which we tacitly apply from now on), we find that 
\begin{align}
I  & \lesssim \big \| \big ( \frac{e^{-\mu |x|}}{|x|} \ast (|u|^2 - |v|^2) \big )(u+v) \big \|_2 + \big \| \big ( \frac{e^{-\mu |x|}}{|x|} \ast ( |u|^2 + |v|^2) \big ) (u-v) \big \|_2 \nonumber \\
& \lesssim \big \|  \frac{e^{-\mu |x|}}{|x|} \ast (|u|^2 - |v|^2) \big \|_6 \| u+v \|_3 + \big \|  \frac{e^{-\mu |x|}}{|x|} \ast ( |u|^2 + |v|^2) \big \|_\infty \| u-v \|_2. \label{ineq-I}
\end{align}  
Observing that $e^{-\mu |x|} |x|^{-1}  \in L^3_w$ holds, the first term of right-hand side of (\ref{ineq-I}) can be bounded by means of the weak Young inequality (see, \eg, \cite{Lieb+Loss2001}) as follows
\begin{equation} \label{ineq-young}
\big \| \frac{e^{-\mu |x|}}{|x|} \ast ( |u|^2 - |v|^2 ) \big \|_6 \lesssim \big \| \frac{e^{- \mu |x|}}{|x|} \big \|_{3,w} \| |u|^2 - |v|^2 \|_{6/5} \lesssim \| u + v \|_3 \| u - v \|_2 .
\end{equation}
The second term occurring in (\ref{ineq-I}) can be estimated by noting that
\begin{align} 
\big \| \frac{e^{-\mu |x|}}{|x|} \ast |u|^2 \big \|_\infty & \lesssim \sup_{y \in \RR^3} \int_{\RR^3} \frac{|u(x)|^2}{|x-y|} \, dx  \lesssim \| (-\DD)^{1/4} u \|_2^2, \label{ineq-kato}
\end{align}
which follows from the operator inequality $|x-y|^{-1} \leq (\pi/2) (-\DD_{x-y})^{1/2}$ (see, \eg, \cite[Section V.5.4]{Kato1980}) and translational invariance, \ie, we use that $\DD_{x-y} = \DD_x$ holds for all $y \in \RR^3$. Combining now (\ref{ineq-young}) and (\ref{ineq-kato}) we find that
\begin{align*}
I & \lesssim \| u + v \|_3^2 \| u - v \|_2 + ( \| u \|_{H^{1/2}}^2 + \| v \|_{H^{1/2}}^2 ) \| u - v\|_2 \nonumber \\ 
& \lesssim ( \| u \|_{H^{1/2}}^2 + \| v \|_{H^{1/2}}^2) \| u - v \|_{H^{1/2}} ,
\end{align*} 
where we make use of the Sobolev inequality $\| u \|_3 \lesssim \| u \|_{H^{1/2}}$ in $\RR^3$. 

It remains to estimate $II$.  To do so, we appeal to the generalized (or fractional) Leibniz rule (see Appendix \ref{app-calc}) leading to
\begin{align}
II & \lesssim \big \| \Ds{1/2}  \big [ \big ( \frac{e^{-\mu |x|}}{|x|} \ast (|u|^2 - |v|^2) \big )(u+v) \big ] \big \|_2  \nonumber \\ 
& \quad + \big \| \Ds{1/2}  \big [ \big ( \frac{e^{-\mu |x|}}{|x|} \ast ( |u|^2 + |v|^2) \big ) (u-v) \big ] \big \|_2  \nonumber \\
& \lesssim \big \| \Ds{1/2} \big ( \frac{e^{-\mu |x|}}{|x|} \ast (|u|^2 - |v|^2) \big ) \big  \|_{6} \| u + v \|_3 \nonumber \\ 
& \quad + \big \|  \frac{e^{-\mu |x|}}{|x|} \ast (|u|^2 - |v|^2) \big  \|_\infty \| \Ds{1/2} (u+v) \|_2 \nonumber \\
& \quad + \big \| \Ds{1/2} \big [ (\frac{e^{-\mu |x|}}{|x|} \ast ( |u|^2 + |v|^2 ) \big ]  \big \|_6 \| u - v \|_3 \nonumber \\
& \quad + \big \| \frac{e^{-\mu |x|}}{|x|} \ast ( |u|^2 + |v|^2 ) \big \|_\infty \| \Ds{1/2} (u-v) \|_2 \label{ineq-II} .
\end{align} 
By referring to Appendix \ref{app-calc}, we notice that $\frac{e^{-\mu |x|}}{4 \pi |x|} \ast f$ can be expressed as $\Ds{-2} f=(\mu^2-\DD)^{-1}f$ in $\RR^3$ (here $f \in \mathcal{S}(\RR^3)$ is initially assumed, but our arguments follow by density). Thus, the first term of the right-hand side of (\ref{ineq-II}) is found to be
\begin{align}
\big \| \Ds{1/2} \big ( \frac{e^{-\mu |x|}}{|x|} \ast (|u|^2 - |v|^2) \big ) \big  \|_{6} & \lesssim \| \Ds{1/2-2} ( |u|^2 - |v|^2 ) \|_6 \nonumber \\
& \lesssim \| \Ds{-3/2} ( |u|^2 - |v|^2 ) \|_6 \nonumber \\
& \lesssim \big \| G^\mu_{3/2} \ast (|u|^2 - |v|^2) \big \|_6 \nonumber \\
& \lesssim \big \| G^\mu_{3/2} \big \|_{2,w} \| |u|^2 - |v|^2 \|_{3/2} \nonumber \\
& \lesssim \| u + v\|_3 \| u - v \|_3, \label{ineq-kato2}
\end{align}
where we use weak Young's inequality together with the fact that $\Ds{-3/2}f$ corresponds to $G_{3/2}^\mu \ast f$ with some $G_{3/2}^\mu \in L^2_w(\RR^3)$; see (\ref{eq-A1}). The $\| \cdot \|_\infty$-part of the second term occurring in (\ref{ineq-II}) can be estimated by using the Cauchy-Schwarz inequality and (\ref{ineq-kato}) once again:
\begin{align}
\big \|  \frac{e^{-\mu |x|}}{|x|} \ast (|u|^2 - |v|^2) \big \|_\infty  & \leq  \big \| \frac{1}{|x|} \ast ( |u|^2 - |v|^2 ) \big \|_\infty \nonumber \\
& \lesssim \sup_{y \in \RR^3} \Big | \int_{\RR^3} \frac{|u(x)|^2 - |v(x)|^2}{|x-y|} \, dx \Big | \nonumber \\
& \lesssim  \sup_{y \in \RR^3} \Big |  \langle (u(x) + v(x)), \frac{1}{|x-y|} (u(x) - v(x)) \rangle \Big | \nonumber \\
& \lesssim  \| (-\DD)^{1/4} (u +v) \|_2 \| (-\DD)^{1/4} (u - v) \|_2 \nonumber \\
& \lesssim  ( \| u \|_{H^{1/2}} + \| v \|_{H^{1/2}}) \| u - v \|_{H^{1/2}} \label{eq-kato2}
\end{align}
The remaining terms in (\ref{ineq-II}) deserve no further comment, since they can be estimated in a similar fashion to all estimates derived so far. Thus, we conclude that
\begin{equation*}
\| J(u) - J(v) \|_{H^{1/2}} \lesssim I + II \lesssim ( \| u \|_{H^{1/2}}^2 + \| v \|_{H^{1/2}}^2 ) \| u - v \|_{H^{1/2}} 
\end{equation*}
and the proof of Lemma \ref{lem-lip} is now complete. \end{proof}

\begin{remarks}
{\em 1) The proof of Lemma \ref{lem-lip} relies on (\ref{ineq-kato}) in a crucial way. Employing just the Sobolev embedding $H^{1/2} \subset L^2 \cap L^3$ (in $\RR^3$) together with the (non weak) Young inequality is not sufficient to conclude that $\| \frac{e^{-\mu |x|}}{|x|} \ast |u|^2 \|_\infty < \infty$ whenever $u \in H^{1/2}$.

2) The proof of Lemma \ref{lem-lip} fails for ``super-critical'' Hartree nonlinearities $J(u) = ( |x|^{-\alpha} \ast |u|^2 ) u$, where $1 < \alpha < 3$. Thus, the choice $\alpha = 1$  represents a borderline case when deriving local Lipschitz continuity in energy space $H^{1/2}$. }
\end{remarks}

\subsubsection*{Step 2: Conclusion}

Returning to the proof of Theorem \ref{th-lwp}, we note that $A$ defined in (\ref{eq-AF}) gives rise to a self-adjoint operator $L^2$ with domain $H^1$. Moreover, its extension to $H^{1/2}$, which we denote by $A : H^{1/2} \rightarrow H^{-1/2}$, generates a $C^0$-group of isometries, $\{ e^{-itA} \}_{t \in \RR}$, acting on $H^{1/2}$. Local well-posedness in the sense of Theorem \ref{th-lwp} now follows by standard methods for evolution equations with locally Lipschitz nonlinearities. That is, existence and uniqueness of a solution $u \in Y_T$ for the integral equation (\ref{eq-integral}) is deduced by a fixed point argument, for $T > 0$ sufficiently small. The equivalence of the integral formulation (\ref{eq-integral}) and the initial value problem (\ref{eq-Hartree}), with $u_0 \in H^{1/2}$, as well as the blow-up alternative can also be deduced by standard arguments; see, \eg, \cite{Pazy1983,Cazenave+Haraux1998} for general theory on semilinear evolution equations. Finally, note that $u \in C^{1}([0,T); H^{-1})$ follows by equation (\ref{eq-Hartree}) itself. The proof of Theorem \ref{th-lwp} is now accomplished.    


\subsection{Proof of Theorem \ref{th-gwp} (Global Well-Posedness)} \label{sec-conservation}

The first step taken in the proof of Theorem \ref{th-gwp} settles conservation of energy and charge that are given by
\begin{equation} \label{eq-defE}
E[u] := \frac{1}{2} \int_{\RR^3}  \overline{u}(x) \Tm u(x) \, dx + \frac{1}{4} \int_{\RR^3} \big ( \frac{\lambda e^{-\mu |x|}}{|x|} \ast |u|^2 \big )(x) \, |u(x)|^2 \, dx,  
\end{equation}
\begin{equation}
N[u] := \int_{\RR^3} |u(x)|^2 \, dx,
\end{equation}
respectively. After deriving the corresponding conservation laws (where proving energy conservation requires a regularization), we discuss how to obtain a-priori bounds on the energy norm of the solution. 


\subsubsection*{Step 1: Conservation Laws}

\begin{lemma} \label{lem-con}
The local-in-time solutions of Theorem \ref{th-lwp} obey conservation of energy and charge, \ie,
\[ E[u(t)] = E[u_0] \quad \mbox{and} \quad N[u(t)] = N[u_0], \]
for all $t \in [0,T)$.
\end{lemma}  

\begin{proof}[Proof of Lemma \ref{lem-con}] Let $u$ be a local-in-time solution derived in Theorem \ref{th-lwp}, and let $T$ be its maximal time of existence. Since $u(t) \in H^{1/2}$ holds, we can multiply (\ref{eq-Hartree}) by $i \bar{u}(t)$ and integrate over $\RR^3$. Taking then real parts yields
\begin{equation}
\frac{d}{dt} N[u(t)] = 0 \quad \mbox{for} \quad t \in [0,T), 
\end{equation}
which shows conservation of charge. 

At a formal level, conservation of energy follows by multiplying (\ref{eq-Hartree}) with $\dot{\bar{u}}(t) \in H^{-1/2}$ and integrating over space, but the paring of two elements of $H^{-1/2}$ is not well-defined. Thus, we have to introduce a regularization procedure as follows; see also, \eg, \cite{Cazenave2003,Ginibre+Velo2000} for other regularization methods for nonlinear (nonrelativistic) Schr\"odinger equations. Let us define the family of operators
\begin{equation}
\Me := (\eps A + 1)^{-1}, \quad \mbox{for} \quad \eps > 0,
\end{equation}  
where the operator $A = \Tm \geq 0$ is taken from (\ref{eq-AF}). Consider the sequences of embedded spaces
\begin{equation}
\ldots  H^{3/2} \hookrightarrow H^{1/2} \hookrightarrow H^{-1/2} \hookrightarrow H^{-3/2}  \ldots
\end{equation}
It is easy to see (by using functional calculus) that the following properties hold.
\begin{enumerate}
\item[a)] For $\eps > 0$ and $s \in \RR$, we have that $\Me$ is a bounded map from $H^{s}$ into $H^{s+1}$.
\item[b)] $\| \Me u \|_{H^{s}} \leq \| u \|_{H^{s}}$ whenever $u \in H^{s}$ and $s \in \RR$.
\item[c)] For $u \in H^{s}$ and $s \in \RR$, we have that $\Me u \rightarrow u$ strongly in $H^{s}$ as $\eps \downarrow 0$.
\end{enumerate} 
We shall use tacitly properties a) -- c) in the following analysis. 

By means of $\Me$ and noting that $E \in C^1(H^{1/2}; \RR)$, we can compute in a well-defined way for $t_1, t_2 \in [0,T)$ as follows
\begin{align}
E[\Me u(t_2)] - E[\Me u(t_1)] & = \int_{t_1}^{t_2} \langle E'(\Me u), \Me \dot{u} \rangle_{H^{-1/2}, H^{1/2}} \, dt \nonumber \\
& = \int_{t_1}^{t_2} \mbox{Re} \, \langle A \Me u + F(\Me u), -i \Me (A u + F(u)) \rangle \, dt \nonumber \\
& = \int_{t_1}^{t_2} \mbox{Im} \, \Big [ \langle A \Me u, \Me A u \rangle + \langle F(\Me u), \Me A u \rangle \nonumber \\
& \quad + \langle A \Me u, \Me F(u) \rangle + \langle F(\Me u), \Me F(u) \rangle \Big ] \, dt \nonumber \\
& =: \int_{t_1}^{t_2} f_\eps(t) \, dt, \label{eq-E}
\end{align}
where we write $u=u(t)$ for brevity and recall the definition of $F$ from (\ref{eq-AF}). We observe that the first term in $f_\eps(t)$ is the ``most singular'' part, \ie, if $\eps = 0$ we would have pairing of two $H^{-1/2}$-elements. But for $\eps > 0$ we can use the obvious fact that $\Me A = A \Me$ holds and conclude that
\begin{equation*}
\mbox{Im} \, \langle A \Me u, \Me A u \rangle = \mbox{Im} \, \langle A \Me u, A \Me u \rangle = 0.
\end{equation*}
Notice that this manipulation is well-defined, because $A \Me u$ and $\Me A u$ are in $H^{1/2}$ whenever $u \in H^{1/2}$. After some simple calculations, we find $f_\eps(t)$ to be of the form
\begin{align*}
f_\eps(t)  = \mbox{Im} \, & \Big [  \langle F(\Me u), \Me A u \rangle + \langle A \Me u, \Me F(u) \rangle + \langle F(\Me u), \Me F(u) \rangle \Big ] \nonumber \\
 = \mbox{Im} \, & \Big [ \langle A^{1/2} F( \Me u), A^{1/2} \Me u \rangle + \langle A^{1/2} \Me u, A^{1/2} \Me F(u) \rangle \nonumber \\
& \quad + \langle F(\Me u), \Me F(u) \rangle \Big ],
\end{align*}
Since $\Me u \rightarrow u$ strongly in $H^{1/2}$ as $\eps \downarrow 0$,  we can infer, by Lemma \ref{lem-lip}, that
\begin{align*}
\lim_{\eps \downarrow 0} f_\eps(t) & = \mbox{Im} \, \big [ \langle A^{1/2} F(u), A^{1/2} u \rangle + \langle A^{1/2} u, A^{1/2} F(u) \rangle + \langle F(u), F(u) \rangle \big ] \nonumber \\
& = \mbox{Im} \, ( \mbox{Real Number} ) = 0.
\end{align*}
To interchange the $\eps$-limit with the $t$-integration in (\ref{eq-E}), we appeal to the dominated convergence theorem. That is, we seek for a uniform bound on $f_\eps(t)$. In fact, by using the Cauchy-Schwarz inequality and Lemma \ref{lem-lip} again we find the following estimate
\begin{align*}
|f_\eps(t)| & \lesssim | \langle A^{1/2} F(\Me u), A^{1/2} \Me u \rangle | + | \langle A^{1/2} \Me u, A^{1/2} \Me F(u) \rangle | \nonumber \\
& \quad + | \langle F(\Me u), \Me F(u) \rangle | \nonumber \\
& \lesssim \| A^{1/2} F(\Me u) \|_2 \| A^{1/2} \Me u \|_2 + \| A^{1/2} \Me u \|_2 \| A^{1/2} \Me F(u) \|_2 \nonumber \\
& \quad + \| F(\Me u) \|_2 \| \Me F(u) \|_2 \nonumber \\
& \lesssim \| u \|^4_{H^{1/2}} + \| u  \|^6_{H^{1/2}} ,
\end{align*} 
for all $\varepsilon > 0$. Putting now all together leads to conservation of energy, \ie, we find for all $t_1, t_2 \in [0,T)$ that
\begin{align*}
E[u(t_2)] - E[u(t_1)] & = \lim_{\eps \downarrow 0} \big ( E[\Me u(t_2)] - E[\Me u(t_1)] \big ) \nonumber \\
& = \lim_{\eps \downarrow 0} \int_{t_1}^{t_2} f_\eps(t) \, dt = \int_{t_1}^{t_2} \lim_{\eps \downarrow 0} f_\eps(t) \, dt = 0 .
\end{align*}
This completes the proof of Lemma \ref{lem-con}. \end{proof}


\subsubsection*{Step 2: A-Priori Bounds}

To fill the last gap towards the global well-posedness result of Theorem \ref{th-gwp}, we now discuss how to obtain a-priori bounds on the energy norm. By the blow-up alternative of Theorem \ref{th-lwp}, global-in-time existence follows from an a-priori bound of the form
\begin{equation} \label{ineq-apriori}
\| u(t) \|_{H^{1/2}} \leq C( u_0 ).
\end{equation}

First, let us assume that $\lambda \geq 0$ holds. Then, for all $t \in [0,T)$, we find from Lemma \ref{lem-con} and (\ref{eq-defE}) that
\begin{equation*}
\| (-\DD)^{1/4} u(t) \|_2 \lesssim E[u(t)] = E[u_0] .
\end{equation*} 
This implies together with charge conservation derived in Lemma \ref{lem-con}, \ie,   
\begin{equation}
\| u(t) \|_2^2 = N[u(t)] = N[u_0] 
\end{equation}
an a-priori estimate (\ref{ineq-apriori}). Therefore condition i) in Theorem \ref{th-gwp} is sufficient for global existence.

Suppose now a focusing nonlinearity, \ie, $\lambda < 0$ holds, and without loss of generality we assume that $\lambda = -1$ is true (the general case follows by rescaling). Now we can estimate as follows.
\begin{align}
E[u] & = \frac{1}{2} \| (-\DD+m^2)^{1/4} u \|_2^2  - \frac{1}{4} \int_{\RR^3} \big ( \frac{e^{-\mu |x|}}{|x|} \ast |u|^2 ) (x) \, |u(x)|^2 \, dx \nonumber \\
& \geq \frac{1}{2} \| (-\DD+m^2)^{1/4} u \|_2^2  - \frac{1}{4} \int_{\RR^3} \big ( \frac{1}{|x|} \ast |u|^2 ) (x) \, |u(x)|^2 \, dx \nonumber \\
& \geq \frac{1}{2} \| (-\DD)^{1/4} u \|_2^2 - \frac{1}{4 K} \| (-\DD)^{1/4} u \|_2^2 \| u \|_2^2 \nonumber \\
& = \left ( \frac{1}{2} - \frac{1}{4 K} \| u \|_2^2 \right ) \| (-\DD)^{1/4} u \|_2^2,
\end{align} 
where $K > 0$ is the best constant taken from Appendix \ref{app-ground}. Thus, energy conservation leads to an a-priori bound on the $H^{1/2}$-norm of the solution, if 
\begin{equation} \label{ineq-K0}
 \| u_0 \|_2^2 < 2 K
\end{equation}
holds. In fact, the constant $K$ satisfies
\[
K = \frac{ \| Q \|_2^2 }{2 }  > \frac{2}{\pi},
\]
where $Q(x)$ is a strictly positive (ground state) solution of 
\begin{equation}
\T Q - ( \frac{1}{|x|} \ast |Q|^2 ) Q = - Q;
\end{equation}      
see Appendix \ref{app-ground}. Going back to (\ref{ineq-K0}), we find that
\begin{equation}
 \| u_0 \|_2^2 < \| Q \|_2^2
\end{equation}
is sufficient for global existence for $\lambda = -1$. The assertion of Theorem \ref{th-gwp} for all $\lambda < 0$ now follows by simple rescaling. The proof of Theorem \ref{th-gwp} is now complete. 

\subsection{Proof of Theorem \ref{th-reg} (Higher Regularity)}

To prove Theorem \ref{th-reg}, we need the following generalization of Lemma \ref{lem-lip}, whose proof is a careful but straightforward generalization of the proof of Lemma \ref{lem-lip}. We defer the details to Appendix \ref{app-lem-lip3}.

\begin{lemma} \label{lem-lip3}
For $\mu \geq 0$ and $s \geq 1/2$, the map $J(u) := ( \frac{e^{-\mu |x|}}{|x|} \ast |u|^2) u$ is locally Lipschitz continuous from $H^s(\RR^3)$ into itself with
\[
\| J(u) - J(v) \|_{H^s} \lesssim ( \| u \|_{H^s}^2 + \| v \|_{H^s}^2 ) \| u - v \|_{H^s}
\]
for all $u, v \in H^s(\RR^3)$. Moreover, we have that 
\[ \| J(u ) \|_{H^s} \lesssim \| u \|^2_{H^r} \| u \|_{H^s} \]
holds for all $u \in H^s(\RR^3)$, where $r = \max \{s-1,1/2 \}$.
\end{lemma}

Local well-posedness of (\ref{eq-HartreeV}) in $H^s$, for $ s > 1/2$, can be shown now as follows. We note that $\{ e^{-itA} \}_{t \in \RR}$, with $A = \Tm$, is a $C^0$-group of isometries on $H^s$. Moreover, since the nonlinearity defined in (\ref{eq-AF}), is locally Lipschitz continuous from $H^s$ into itself, local well-posedness in $H^s$ follows similarly as explained in the proof of Theorem \ref{th-lwp} for $H^{1/2}$. To show global well-posedness in $H^s$, we prove by induction and Lemma \ref{lem-lip3} that an a-priori bound on the $H^{1/2}$-norm of solution implies uniform bounds on the $H^s$-norm on any compact interval $[0,T_*] \subset [0,T)$. This claim follows from (\ref{eq-integral}) and the second inequality stated in Lemma \ref{lem-lip3} by noting that
\begin{align*}
\| u(t) \|_{H^s} & \leq \| e^{-itA} u_0 \|_{H^s} + \int_0^t \| e^{-i(t-\tau) A} F(u(\tau)) \|_{H^{s}} \, d\tau \\
& \leq \| u_0 \|_{H^s} + \int_0^t \| F(u(\tau)) \|_{H^s} \, d \tau \\
& \lesssim C_1 + C_2 \int_0^t \| u(\tau) \|_{H^s} \, d\tau, 
\end{align*}  
holds, provided that $\|u(t)\|_{H^r} \lesssim 1$ for $r = \max \{s-1,1/2 \}  < s$. Invoking Gronwall's inequality we conclude that
\[ \| u(t) \|_{H^s} \lesssim  e^{C_2 T_*}, \quad \mbox{for} \quad t \in [0,T_*] \subset [0,T). \]
Induction now implies that an a-priori bound on $\| u(t) \|_{H^{1/2}}$ guarantees uniform bounds $\| u(t) \|_{H^s}$ on any compact interval $I \subset [0,T)$. Thus, the maximal time of existence of an $H^s$-valued solution coincides with the maximal time of existence when viewed as an $H^{1/2}$-valued solution. Therefore sufficient conditions for global existence for $H^{1/2}$-valued solutions imply global-in-time $H^s$-valued solutions. This completes the proof of Theorem \ref{th-reg}.

\subsection{Proof of Theorem \ref{th-V} (External Potentials)} \label{sec-proof-V}

Let $V=V_+ + V_-$ satisfy Assumption \ref{ass-V} in Section \ref{sec-main}. We introduce the quadratic form
\begin{equation} \label{eq-quad}
\mathcal{Q}(u,v) := \langle u, \Tm v \rangle + \langle u, V_- v \rangle + \langle u, V_+ v \rangle, 
\end{equation}
which is well-defined on the set (energy space)
\begin{equation} \label{eq-X}
X:= \big \{ u \in L^2(\RR^3) : \mathcal{Q}(u,u) < \infty \big \} .
\end{equation}
Note that Assumption \ref{ass-V} also guarantees that $C_0^\infty(\RR^3) \subset X$. It easy to show that our assumption on $V$ implies that the quadratic form (\ref{eq-quad}) is bounded from below, \ie, we have $\mathcal{Q}(u,u) \geq -M \langle u, u \rangle$ holds for all $u \in X$ and some constant $M \geq 0$. By the semi-boundedness of $\mathcal{Q}$, we can assume from now on (and without loss of generality) that 
\begin{equation}
\mathcal{Q}(u,u) \geq 0
\end{equation}
holds for all $u \in X$. Since $\mathcal{Q}(\cdot, \cdot)$ is closed (it is a sum of closed forms), the energy space $X$ equipped with its norm
\begin{equation}
\| u \|_X := \sqrt{\langle u, u \rangle + \mathcal{Q}(u,u)}
\end{equation} 
is complete, and we have the equivalence 
\begin{equation} \label{eq-eqX}
\| u \|_{H^{1/2}} + \| V_+^{1/2} u \|_2 \lesssim \| u \|_X \lesssim \| u \|_{H^{1/2}} + \| V^{1/2}_+ u \|_2 .
\end{equation}
Furthermore, there exists a nonnegative, self-adjoint operator
\begin{equation}
A : D(A) \subset L^2 \rightarrow L^2
\end{equation}
with $X=D(A^{1/2})$, such that 
\begin{equation}
\langle u, A v \rangle = \mathcal{Q}(u,v)
\end{equation} 
holds for all $u \in X$ and $v \in D(A)$; see, \eg, \cite{Kato1980}. This operator can be extended to a bounded operator, still denoted by $A : X \rightarrow X^*$, where $X^*$ is the dual space of $X$. 

To prove now the assertion about local well-posedness in Theorem \ref{th-V}, we have to generalize Lemma \ref{lem-lip} to the following statement.

\begin{lemma} \label{lem-lip2}
Suppose $\mu \geq 0$ and let $V$ satisfy Assumption \ref{ass-V}. Then the map $J(u) := ( \frac{e^{-\mu |x|}}{|x|} \ast |u|^2) u$ is locally Lipschitz continuous from $X$ into itself with
\[
\| J(u) - J(v) \|_X \lesssim ( \| u \|_X^2 + \| v \|_X^2 ) \| u - v \|_X
\]
for all $u,v \in X$.
\end{lemma}

\begin{proof}[Proof of Lemma \ref{lem-lip2}]
By (\ref{eq-eqX}), it suffices to estimate $\| J (u) - J(v) \|_{H^{1/2}}$ and $\| V^{1/2}_+ [ J(u) - J(v) ] \|_2$ separately. By Lemma \ref{lem-lip}, we know that
\begin{align*}
\| J(u) - J(v) \|_{H^{1/2}} & \lesssim ( \| u \|_{H^{1/2}}^2 + \| v \|_{H^{1/2}}^2 ) \| u - v \|_{H^{1/2}} \\
& \lesssim  ( \| u \|_{X}^2 + \| v \|_{X}^2 ) \| u - v \|_{X} .
\end{align*}
It remains to estimate $\| V^{1/2}_+ [J(u)-J(v)] \|_2$, which can be achieved by recalling (\ref{eq-kato2}) and proceeding as follows.
\begin{align}
\| V^{1/2}_+ [J(u)-J(v)] \|_2 & \lesssim \big \| V^{1/2}_+ \big [ \big (\frac{e^{-\mu |x|}}{x} \ast ( |u|^2 - |v|^2) \big) (u+v) \big ] \big \|_2 \nonumber \\
& \quad + \big \| V^{1/2}_+ \big [ \big (\frac{e^{-\mu |x|}}{x} \ast ( |u|^2 + |v|^2) \big) (u-v) \big ] \big \|_2 \nonumber \\
& \lesssim \big \| \frac{e^{-\mu |x|}}{|x|} \ast ( |u|^2 - |v|^2) \big \|_\infty \| V^{1/2}_+ (u+v) \|_2 \nonumber  \\
& \quad + \big \| \frac{e^{-\mu |x|}}{|x|} \ast ( |u|^2 + |v|^2 ) \big \|_\infty \| V^{1/2}_+ (u-v) \|_2 \nonumber \\
& \lesssim \big \| u + v \|_{H^{1/2}} \| u -v \|_{H^{1/2}} \| V^{1/2}_+ (u+v) \|_2 \nonumber \\
& \quad + ( \| u \|^2_{H^{1/2}} + \| v \|^2_{H^{1/2}} ) \| V^{1/2}_+ (u-v) \|_2 \nonumber \\
& \lesssim \big ( \| u \|^2_X + \| v \|^2_X ) \| u -v \|_X \nonumber .
\end{align}   
This completes the proof of Lemma \ref{lem-lip2}. \end{proof}
 
Returning to the proof of Theorem \ref{th-V}, we simply note that $\{ e^{-it A} \}_{t \in \RR}$ is a $C^0$-group of isometries on $X$, where $A = \Tm + V$ is defined in the form sense (see above). By Lemma \ref{lem-lip2}, the nonlinearity is locally Lipschitz on $X$. Thus, local well-posedness now follows in the same way as for Theorem \ref{th-lwp}.

To establish global well-posedness we have to prove conservation of charge, $N[u]$, and energy, $E[u]$, which is for (\ref{eq-HartreeV}) given by
\begin{align}
E[u]  & :=   \frac{1}{2} \int_{\RR^3}  \overline{u}(x) \Tm u(x) \, dx + \frac{1}{2} \int_{\RR^3} V(x) | u(x) |^2 \, dx  \nonumber \\
& \quad + \frac{1}{4} \int_{\RR^3} \big ( \frac{\lambda e^{-\mu |x|}}{|x|} \ast |u|^2 \big )(x) \, |u(x)|^2  \, dx .
\end{align}
As done in Section \ref{sec-conservation}, we have to employ a regularization method using the class of operators
\begin{equation}
\Me := (\eps A + 1)^{-1}, \quad \mbox{for} \quad \eps > 0,
\end{equation} 
where we assume without loss of generality that $A \geq 0$ holds. The mapping $\Me$ acts on the sequence of embedded spaces
\begin{equation}
\ldots X^{+2} \hookrightarrow X^{+1} \hookrightarrow X^{-1} \hookrightarrow X^{-2} \ldots,
\end{equation}
with corresponding norms given by $\| u \|_{X^s} := \| (1+A)^{s/2} u \|_2$. Note that $X = X^{+1}$ (with equivalent norms) and that its dual space obeys $X^* = X^{-1}$. By using functional calculus, it is easy to show that $\Me$ exhibits properties that are analog to a) -- c) in Section \ref{sec-conservation}. 

The rest of the argument for proving conservation of energy carries over from Section \ref{sec-conservation} without major modifications. Finally, we mention that deriving a-priori bounds on $\| u(t) \|_X$ leads to a similar discussion as presented in Section \ref{sec-conservation}, while noting that we have to take care that $V_-$ has a relative $(-\DD)^{1/2}$-form bound, $0 \leq a < 1$,  introduced in Assumption \ref{ass-V}. This completes the proof of Theorem \ref{th-V}. 

\subsubsection*{Acknowledgments} 

The author is grateful to Demetrios Christodoulou, J\"urg Fr\"ohlich, Lars Jonsson, and Simon Schwarz for many valuable and inspiring discussions.

\begin{appendix}

\section{Appendix} \label{sec-app}

\subsection{Fractional Calculus} \label{app-calc}

The following result (generalized Leibniz rule) is proved in \cite{Gulisashvili+Kon1996} for Riesz and Bessel potentials of order $s \in \RR$, which are denoted by $(-\DD)^{s/2}$ and $(1-\DD)^{s/2}$, respectively. But as a direct consequence of the Milhin multiplier theorem \cite{Bergh+Loefstroem1976}, the cited result holds for $\mathcal{D}^s := (\mu^2 - \DD)^{s/2}$, where $\mu \geq 0$ is a fixed constant.   
 
\begin{lemma*}[Generalized Leibniz Rule]  Suppose that $1 < p < \infty$, $s \geq 0$, $\alpha \geq 0$, $\beta \geq 0$, and $1/p_i + 1/q_i = 1/p$ with $i=1,2$, $1 < q_i \leq \infty$, $1 < p_i \leq \infty$. Then
\[ \| \mathcal{D}^s (fg) \|_p \leq c ( \| \mathcal{D}^{s+\alpha} f \|_{p_1} \|  \mathcal{D}^{-\alpha} g \|_{q_1} + \| \mathcal{D}^{-\beta} f \|_{p_2} \| \mathcal{D}^{s + \beta} g \|_{q_2} ) , \]
where the constant $c$ depends on all of the parameters but not on $f$ and $g$. 
\end{lemma*}

A second fact we use in the proof of our main result is as follows. For $0 < \alpha < 3$ and $\mu \geq 0$, the potential operator $\mathcal{D}^{-\alpha}=(\mu^2-\DD)^{-\alpha/2}$ corresponds to $f \mapsto G_\alpha^{\mu} \ast f$, with $f \in \mathcal{S}(\RR^3)$, and we have that
\begin{equation} \label{eq-A1}
G_\alpha^{\mu} \in L^{3/(3-\alpha)}_w(\RR^3). 
\end{equation}   
To see this, we refer to the inequality and the exact formula
\begin{equation}
0 \leq  G_\alpha^\mu(x) \leq G_\alpha^{0}(x) = \frac{c_\alpha}{|x|^{3-\alpha}}, \quad \mbox{for} \quad \mu \geq 0 \quad \mbox{and} \quad  0 < \alpha < 3,
\end{equation}
with some constant $c_\alpha$; these facts can be derived from \cite[Section V.3.1]{Stein1970}. Now (\ref{eq-A1}) follows from $|x|^{-\sigma} \in L^{3/\sigma}_w(\RR^3)$ whenever $0 < \sigma < 3$. Another observation used in Section 2 is the well-known explicit formula
\begin{equation}
G^\mu_2(x) = \frac{e^{-\mu |x|}}{4 \pi |x|}.
\end{equation}
That is, $(\mu^2-\DD)$ in $\RR^3$ has the Green's function $\frac{e^{-\mu |x|}}{4 \pi |x|}$ with vanishing boundary conditions.

\subsection{Ground States} \label{app-ground}

We consider the functional (see also \cite{Lieb+Yau1987})
\begin{equation}
K[u] := \frac{ \| (-\DD)^{1/4} u \|_2^2 \| u \|_2^2}{ \int_{\RR^3} ( |x|^{-1} \ast |u|^2 )(x) \, |u(x)|^2 \, dx},
\end{equation}
which is well-defined for all $u \in H^{1/2}$ with $u \not \equiv 0$. Note that by using (\ref{ineq-kato}) we can estimate the denominator in $K[u]$ as follows.
\begin{equation} \label{ineq-K0_}
\int_{\RR^3} \big ( \frac{1}{|x|} \ast |u|^2) (x) \, |u(x)|^2 \, dx \leq \big \| \frac{1}{|x|} \ast |u|^2 \big \|_\infty \| u \|_2^2 \leq \frac{\pi}{2} \| (-\DD)^{1/4} u \|_2^2 \| u \|_2^2,
\end{equation}
which leads to the bound
\begin{equation}
\frac{2}{\pi} \leq K[u] < \infty .
\end{equation}  
Indeed, we will see that the estimate from below is a strict inequality. With respect to the related variational problem
\begin{equation} \label{eq-Kinf}
K := \inf \big \{ K[u] : u \in H^{1/2}(\RR^3), \; u \not \equiv 0 \big  \} 
\end{equation} 
we can state the following result.

\begin{lemma*}[Ground States]
There exists a minimizer, $Q \in H^{1/2}(\RR^3)$, for (\ref{eq-Kinf}), and we have the following properties.
\begin{enumerate}
\item[i)] $Q(x)$ is a smooth function that can be chosen to be real-valued, strictly positive, and spherically symmetric with respect to the origin. It satisfies
\begin{equation} \label{eq-Q3}
 \T Q - \big ( \frac{1}{|x|} \ast |Q|^2) Q = - Q,
\end{equation}
and it is nonincreasing, \ie, we have that $Q(x) \geq Q(y)$ whenever $|x| \leq |y|$.
\item[ii)] The infimum satisfies $K = \| Q \|_2^2/2$ and $K > 2/\pi$.   
\end{enumerate}
\end{lemma*}

\begin{proof}[Sketch of Proof]
We present the main ideas for the proof of the preceding lemma. That (\ref{eq-Kinf}) is attained at some real-valued, radial, nonnegative and nonincreasing function $Q(x) \geq 0$ can be proved by direct methods of variational calculus and rearrangement inequalities; see also \cite{Weinstein1983} for a similar variational problem for nonrelativistic Schr\"odinger equations with local nonlinearities. Furthermore, any minimizer, $Q \in H^{1/2}$, has to satisfy the corresponding Euler-Lagrange equation that reads
\begin{equation} \label{eq-Q}
\T Q - \big ( \frac{\lambda}{|x|} \ast |Q|^2) Q = - Q,
\end{equation}
after a suitable rescaling $Q(x) \mapsto a Q( b x)$ with some $a, b > 0$. 

Let us make some comments about the properties of $Q$. Using an bootstrap argument and Lemma \ref{lem-lip3} for the nonlinearity, it follows that $Q$ belongs to $H^s$, for all $s \geq 1/2$. Hence it is a smooth function. To see that $Q(x) \geq 0$ is strictly positive, \ie, $Q(x) > 0$, we rewrite equation (\ref{eq-Q}) such that
\begin{equation} \label{eq-Q2}
Q = \big  ( \sqrt{-\DD} + 1 \big )^{-1} W,
\end{equation}
where $W := ( |x|^{-1} \ast |Q|^2 ) Q$. By functional calculus, we have that
\begin{equation} \label{eq-Res}
\big (\sqrt{-\DD} + 1 \big )^{-1} = \int_{0}^\infty e^{-t} e^{-t \sqrt{-\DD}} \, dt.
\end{equation} 
Next, we notice by the explicit formula for the kernel (in $\RR^3$) 
\[
e^{-t \sqrt{-\DD}}(x,y) = \mathcal{F}^{-1} \big (e^{-t |\xi|} \big )(x-y) =  C \cdot \frac{t}{[t^2+ |x-y|^2]^2},
\]
with some contant $C > 0$; see, \eg, \cite{Lieb+Loss2001}. This explicit formula shows that $e^{-t \sqrt{-\DD}}$ is positivity improving. This means that if $f \geq 0$ with $f \not \equiv 0$ then $e^{-t\sqrt{-\DD}} f >0$ almost everywhere. Hence $(\sqrt{-\DD} + 1)^{-1}$ is also positivity improving, by (\ref{eq-Res}), and we conclude that $Q(x) > 0$ holds almost everywhere, thanks to (\ref{eq-Q2}) and $W \geq 0$. Moreover, we know that $Q(x)$ is a nonincreasing, continuous function. Therefore $Q(x) > 0$ holds in the strong sense, \ie, for every $x \in \RR^3$.  

Finally, to see that ii) holds, we consider the variational problem 
\begin{equation} \label{eq-EQ}
I_N  := \inf \big \{  E[u] : u \in H^{1/2}(\RR^3), \; \| u \|_2^2 = N \big \}, 
\end{equation}
where $N > 0$ is a given parameter and
\[
E[u] = \frac{1}{2} \| (-\DD)^{1/4} u \|^2_2 - \frac{1}{4} \int_{\RR^3} \big ( \frac{1}{|x|} \ast |u|^2) \, |u(x)|^2 \, dx .
\]
Due to the scaling behavior $E[\alpha^{3/2} u( \alpha \cdot )] = \alpha E[u]$, we have that either $I_N = 0$ or $I_N = -\infty$ holds. By noting that
\begin{equation*}
E[u] \geq \big ( \frac{1}{2} - \frac{N}{4K} \big ) \| (-\DD)^{1/4} u \|_2^2 , 
\end{equation*}
and the fact that equality holds if and only if $u$ minimizes $K[u]$, we find that $I_N = 0$ holds if and only if $N \leq N_c := 2K$. Moreover, $I_N=0$ is attained if and only if $N = N_c$. Let $\tilde{Q}$ be such a minimizer with $\| \tilde{Q} \|_2^2 = N_c$. Thanks to the proof of part i), we can assume without loss of generality that $\tilde{Q}$ is real-valued, radial, and strictly positive. Calculating the Euler-Lagrange equation for (\ref{eq-EQ}), with $N=N_c$, yields
\[
\T \tilde{Q} - \big ( \frac{1}{|x|} \ast |\tilde{Q}|^2 \big ) \tilde{Q} = -\theta \tilde{Q},
\]
for some multiplier $\theta$, where it is easy to show that $\theta > 0$ holds. Putting now $Q(x) = \theta^{-3/2} \tilde{Q}(\theta^{-1} x)$, which conserves the $L^2$-norm, leads to a ground state $Q(x)$ satisfying (\ref{eq-Q}). Thus, we have that 
\[ K = \| \tilde{Q} \|_2^2 / 2 = \| Q \|_2^2 /2 . \] 

To prove that $K > 2/\pi$ holds, let us assume $K = \pi /2$. This implies that the first inequality in (\ref{ineq-K0_}) is an equality for $u(x) = Q(x) > 0$. But this leads to $(|x|^{-1} \ast |Q|^2)(x)=\mbox{const.}$, which is impossible. \end{proof}

\subsection{Proof of Lemma \ref{lem-lip3}} \label{app-lem-lip3}

\begin{proof}[Proof of Lemma \ref{lem-lip3}]
We only show the second inequality derived in Lemma \ref{lem-lip3}, since the first one can be proved in a similar way. 

Let $\mu \geq 0$ and $s \geq 1/2$. We put $\Ds{\alpha} := (\mu^2 - \DD)^{\alpha/2}$ for $\alpha \in \RR$. By the generalized Leibniz rule and (\ref{ineq-kato}), we have that
\begin{align}
 \| \Ds{s} J(u) \|_2 & \lesssim \| \Ds{s} [ (\Ds{-2} |u|^2) u ] \|_2 \nonumber \\
 & \lesssim \| \Ds{s-2} |u|^2 \|_{p_1} \| u \|_{q_1}  + \| \Ds{-2} |u|^2 \|_\infty \| \Ds{s} u \|_2  \nonumber \\
& \lesssim \| \Ds{s-2} |u|^2 \|_{p_1} \| u \|_{q_1} + \| u \|^2_{H^{1/2}} \| u \|_{H^s}, \label{ineq-Hs}
\end{align}
where $1/p_1 + 1/q_1 = 1/2$ with $1 < p_1, q_1 \leq \infty$. The first term of the right-hand side of (\ref{ineq-Hs}) can be controlled as follows, where we introduce $ r = \max \{s-1,1/2 \}$.
\begin{enumerate}
\item[i)] For $1/2 \leq s < 3/2$, we choose $p_1 = 3/s$ and $q_1 = 6/(3-2 s)$ which leads to
\begin{align*}
 \| \Ds{s-2} |u|^2 \|_{3/s} \| u \|_{6/(3-2s)}  & \lesssim \| G^\mu_{2-s} \|_{3/(1+s),w} \| |u|^2 \|_{3/2} \| u \|_{H^s}  \\
 & \lesssim \| u \|_{H^{1/2}}^2 \| u \|_{H^s} \lesssim \| u \|^2_{H^r} \| u \|_{H^s} , 
\end{align*}
where we use the weak Young inequality, Sobolev's inequality $\| u \|_{6/(3-2s)} \lesssim \| u \|_{H^s}$ in $\RR^3$, and (\ref{eq-A1}) once again.
\item[ii)] For $s \geq 3/2$, we put $p_1 = 6$ and $q_1 = 3$ and find 
\begin{align*}
   \| \Ds{s-2} |u|^2 \|_6 \| u \|_{3} & \lesssim \| \Ds{s-3/2} |u|^2 \|_2 \| u \|_3 \lesssim \| \Ds{s-3/2} u \|_6 \| u \|_3^2 \\
   & \lesssim \| \Ds{s-1} u \|_2 \| u \|_3^2 \lesssim \| u \|^2_{H^{r}} \| u \|_{H^s},
\end{align*}
while using twice Sobolev's inequality $\| f \|_6 \lesssim \| \Ds{1/2} f \|_2$ in $\RR^3$.
\end{enumerate}
Putting now all together, we conclude that
\begin{align*} 
\| J(u) \|_{H^s}  & \lesssim \| J(u) \|_2 + \| \Ds{s} J(u) \|_2 \\
& \lesssim \| u \|_{H^{1/2}}^2 \| u \|_2 + \| u \|_{H^r}^2 \| u \|_{H^s}  \lesssim \| u \|_{H^r}^2 \| u \|_{H^s}. 
\end{align*} 
\end{proof}
\end{appendix}

\bibliographystyle{amsalpha}
\bibliography{GWPBib}


\bigskip
\noindent
{\sc Enno Lenzmann \\ Department of Mathematics \\ ETH Z\"urich-Zentrum, HG G 33.1, R\"amistrasse 101 \\ CH-8092 Z\"urich, Switzerland}\\
{\em E-mail:} {\tt lenzmann@math.ethz.ch}

\end{document}